\documentclass[11pt]{article}
\usepackage{amsthm,amssymb,amsmath, color}
\usepackage{amsfonts}
\usepackage{CJK}
\usepackage{mathrsfs}
\usepackage{amsfonts,amssymb}
\usepackage{amsmath,amscd,mathptm}
\usepackage[all]{xy}
\usepackage{pstricks}
\usepackage{fancyhdr}
\usepackage{wasysym, stmaryrd}
\usepackage[dvipdfm,
            colorlinks=ture,
            citecolor=blue,
            linkcolor=black]{hyperref}
\begin{document}
\sloppy
\newcommand{\dickebox}{{\vrule height5pt width5pt depth0pt}}
\newtheorem{Def}{Definition}[section]
\newtheorem{Bsp}{Example}[section]
\newtheorem{Prop}[Def]{Proposition}
\newtheorem{Theo}[Def]{Theorem}
\newtheorem{Lem}[Def]{Lemma}
\newtheorem{Koro}[Def]{Corollary}
\newtheorem{Rem}[Def] {Remark}
\newcommand{\lra}{\longrightarrow}
\newcommand{\ra}{\rightarrow}
\newcommand{\F}{\mathcal {F}}
\newcommand{\Hom}{{\rm Hom}}
\newcommand{\End}{{\rm End}}
\newcommand{\Ext}{{\rm Ext}}
\newcommand{\Tor}{{\rm Tor}}
\newcommand{\pd}{{\rm proj.dim}}
\newcommand{\inj}{{\rm inj}}
\newcommand{\lgd}{{l.{\rm gl.dim}}}
\newcommand{\gld}{{\rm gl.dim}}
\newcommand{\fd}{{\rm fin.dim}}
\newcommand{\Fd}{{\rm Fin.dim}}
\newcommand{\lfd}{l.{\rm Fin.dim}}
\newcommand{\rfd}{r.{\rm Fin.dim}}
\newcommand{\Mod}{{\rm Mod}}
\newcommand{\Proj}{{\rm Proj}}
\newcommand{\modcat}[1]{#1\mbox{{\rm -mod}}}
\newcommand{\pmodcat}[1]{#1\mbox{{\rm -proj}}}
\newcommand{\Pmodcat}[1]{#1\mbox{{\rm -Proj}}}
\newcommand{\injmodcat}[1]{#1\mbox{{\rm -inj}}}
\newcommand{\E}{{\rm E}_{\mathcal {F}}^{{\rm F},\Phi}}
\newcommand{\X}{ \mathscr{X}_{\mathcal {F}}^{{\rm F},\Phi}}
\newcommand{\Y}{\rm \mathscr{Y}_{\mathcal {F}}^{{\rm F},\Phi}}
\newcommand{\A}{\mathcal {A}}
\newcommand{\C}{\rm \mathscr{C}}
\newcommand{\K}{\rm \mathscr{K}}
\newcommand{\D}{\rm \mathscr{D}}
\newcommand{\opp}{^{\rm op}}
\newcommand{\otimesL}{\otimes^{\rm\bf L}}
\newcommand{\otimesP}{\otimes^{\bullet}}
\newcommand{\rHom}{{\rm\bf R}{\rm Hom}}
\newcommand{\projdim}{\pd}
\newcommand{\stmodcat}[1]{#1\mbox{{\rm -{\underline{mod}}}}}
\newcommand{\Modcat}[1]{#1\mbox{{\rm -Mod}}}
\newcommand{\modcatr}[1]{\mbox{{\rm mod}#1}}
\newcommand{\Modcatr}[1]{\mbox{{\rm Mod}#1}}
\newcommand{\Pmodcatr}[1]{\mbox{{\rm Proj}#1}}
\newcommand{\procat}[1]{#1\mbox{{\rm -proj}}}
\newcommand{\Tr}{\rm Tr}
\newcommand{\add}{{\rm add}}
\newcommand{\Imf}{{\rm Im}}
\newcommand{\Ker}{{\rm Ker}}
\newcommand{\EA}{{\rm E^\Phi_\mathcal {A}}}
\newcommand{\pro}{{\rm pro}}
\newcommand{\Coker}{{\rm Coker}}
\newcommand{\id}{{\rm id}}
\renewcommand{\labelenumi}{\Alph{enumi}}
\newcommand{\M}{\mathcal {M}}
\newcommand{\Mf}{\rm \mathcal {M}^f}
\newcommand{\rad}{{\rm rad}}
\newcommand{\injdim}{{\rm inj.dim}}
{\large \bf
\begin{center}Derived equivalences between matrix subrings and their
applications
\end{center}}
\medskip
\centerline{\bf Yiping Chen }
\medskip \centerline{School of Mathematical Sciences, Beijing Normal University}
\medskip
\centerline{100875 Beijing, People's Republic of China}
\medskip
\centerline{E-mail: ypchen@mail.bnu.edu.cn}

\renewcommand{\thefootnote}{\alph{footnote}}
\setcounter{footnote}{-1} \footnote{2000 Mathematics Subject
Classification: 18G20, 18E30; 18G05}

\setcounter{footnote}{-1}\footnote{Keywords: derived category;
derived equivalence; finitistic dimension; global dimension.}

\abstract{In this paper, we construct derived equivalences between
matrix subrings. As applications, we calculate the global dimensions
and the finitistic dimensions of some matrix subrings. And we show
that the finitistic dimension conjecture holds for a class of Harada
algebras and a class of tiled triangular algebras. }

\section{Introduction}

Derived equivalences preserve many homological properties of
algebras such as the number of simple modules, the finiteness of
global dimension and finitistic dimension, the algebraic K-theory
and Hochschild (co)homological groups (see \cite{MB, DS, K, R1, R2,
PX}). Thus, in order to study some homological properties of a given
algebra, we can turn to the one which is derived equivalent to it.

Recently, Hu and Xi have exhibited derived equivalent endomorphism
rings induced by $\mathcal{D}$-split sequences. We find that
$\mathcal{D}$-split sequences give a way to construct derived
equivalences between matrix subrings. In this paper, we will study
the derived equivalences having a characteristic that one of two
rings has relatively simple structure.

As applications, we first investigate the global dimension of a
matrix subring. By the definition of global dimension, Kirkam and
Kuzmanovich in \cite{KK} have calculated the global dimensions of
some matrix subrings. Cowley extended some of their results by
triangular decomposition \cite{C}. As never before, we investigate
some cases by the method of derived equivalences.

Second, we study the finitistic dimension of a matrix subring. For a
ring $A$, the finitistic dimensions are defined as follows:
$\lfd(A)$ is the supremum of the projective dimensions of left
$A$-modules of finite projective dimensions, and $\fd(A)$ is the
supremum of the projective dimensions of finitely generated left
$A$-modules of finite projective dimensions. Kirkman and Kuzmanovich
compute $\lfd(A)$ for a noncommutative Noetherian ring $A$ in
\cite{KK}. By derived equivalences, we calculate $\lfd(A)$ for a
matrix subring $A$. This result is helpful to study the finitistic
dimension conjecture which states that for an Artin algebra $A$,
$\fd(A)$ is finite. This conjecture is still open. We refer the
reader to \cite{XX} on some new advances on this conjecture.

Little is known about whether the finitistic dimension conjecture
holds for matrix subalgebras. Note that the Artin algebra $A$ and
the matrix algebra $M_n(A)$ are Morita equivalent. Thus, in order to
prove that $\fd(A)$ is finite, it is equivalent to prove that
$\fd(M_n(A))$ is finite. Our ideal in this direction is to
investigate the finitistic dimension of a matrix subalgebra. If the
finitistic dimension of $A$ is finite, what could we say about the
finitistic dimension of a matrix subalgebra?

In order to describe the main result precisely, we fix some
notations.

Let $A$ be a ring with identity. Let $A_{i} (2 \leq i \leq n)$ be a
family of subrings of $A$ with the same identity with $A$, and let
$I_{i},\, I_{ij},\, 2\leq i\leq n,\, 2\leq j\leq n-1$ be ideals of
$A$ satisfying that $I_n\subseteq I_{n-1}\subseteq\cdots\subseteq
I_2,\, I_i\subseteq A_i, I_j\subseteq I_{ij},\,
\displaystyle{\sum_{l=j+1}^{i-1}I_{il}I_{lj}}\subseteq I_{ij},
\,2\leq i\leq n, \,2\leq j\leq n-1$. In this way, we can construct
two rings
$$\Lambda =\footnotesize
 \begin{pmatrix}
A & I_2 &I_3 & \cdots & I_{n-1} & I_n\\
A & A_2 & I_3 & \cdots & I_{n-1} & I_n\\
A & I_{32}& A_3& \cdots & I_{n-1} & I_n\\
A & I_{42} & I_{43}& A_4 & \cdots & I_n\\
\vdots & \vdots &\ddots &\ddots  & \ddots& \vdots\\
A & I_{n2} & \cdots & \cdots & I_{n,n-1} & A_n
 \end{pmatrix}\mbox{ and }
\Sigma =\footnotesize
\begin{pmatrix}
A_2/I_2 & 0 &  & \cdots &  &0\\
I_{32}/I_2 &A_3/I_3 & \ddots& & &\vdots\\
I_{42}/I_2 &I_{43}/I_3& A_4/I_4 & & &\\
\vdots & \vdots & \vdots & \ddots& & \\
I_{n2}/I_2 & I_{n3}/I_{3} & I_{n4}/I_4 & \cdots & A_n/I_n &0\\
A/I_2 & A/I_3 & A/I_4 & \cdots & A/I_n & A
\end{pmatrix}$$ with identities. Throughout this paper, $\Lambda$ and $\Sigma$ are rings of this forms.

The main result in this paper is the following:

\begin{Theo}\label{Theo}The two rings $\Lambda$ and $\Sigma$ are derived equivalent.
\end{Theo}

As a direct consequence of Theorem 1.1, we have the following
corollary.

%
%

\begin{Koro}
$(1)$ Let $\Lambda$ be as in Theorem \ref{Theo}. Then
$$\lfd (A)-1\leq \lfd(\Lambda)\leq n+\sum\limits_{i=2}^n\lfd
(A_i/I_i)+\lfd (A).$$

$(2)$ Let $\Lambda$ as in Theorem \ref{Theo}. Then
$$\begin{array}{ll} max\{\lgd (A_i/I_i), \lgd (A), 2\leq i\leq n\}-1\leq
\lgd(\Lambda) \leq \displaystyle{\sum_{i=2,3,\cdots,n} \lgd
(A/I_i)}\\+\lgd (A)+n.
\end{array}$$
\end{Koro}

We define a class of algebras which are called general block
extensions of rings with respect to the decomposition of identity. And we
calculate their global dimensions and finitistic dimensions. We also
get a class of Harada algebras and a class of tiled triangular rings
which satisfy the finitistic dimension conjecture.

This paper is arranged as follows. In section 2, we fix some
 notations and recall some definitions and lemmas needed in this paper.
 In section 3, the proof of the main result is given. In section 4, we give some
 applications of the main result. The definition of general block extensions of
 rings with respect to the decomposition of identity is proposed. We calculate their global dimensions and finitistic
 dimensions. And we also get some classes
 of algebras which satisfy the finitistic dimension conjecture.
 In section 5, we display some examples to illustrate the applications.

\section{Preliminaries}

In this section, we shall recall some basic definitions and results
needed in this paper.

Let $A$ be a ring with identity. We denote $A$-Mod the category of
left $A$-modules and by $A$-mod the category of all finitely
generated left $A$-modules. Mod-$A$ means the category of right
$A$-modules. Given an $A$-module $M$, we denote by $\pd(M)$ the
projective dimension of $M$. The left global dimension of $M$,
denoted by $\lgd(A)$, is the supremum of all $\pd (M)$ with $M\in$
$A$-Mod. By $\add(M)$, we shall mean the full subcategory of
$A$-Mod, whose objects are summands of finite sums of $M$. For two
morphisms $f:{}X \ra Y$ and $g:{}Y \ra Z$ in $A$-Mod, the
composition of $f$ and $g$ is written as $fg$, which is a morphism
form $X$ to $Z$.


Let $A$ be an Artin algebra. A {\em complex} $X^\bullet=(X^i,
d_X^i)$ of $A$-modules is a sequence of $A$-modules and $A$-module
homomorphisms $d^i_X: X^i\ra X^{i+1}$ such that $d^i_Xd^{i+1}_X=0$
for all $i\in \mathbb{Z}$. A {\em morphism} $f^\bullet: X^\bullet\ra
Y^\bullet$ between two complexes $X^\bullet$ and $Y^\bullet$ is a
collection of homomorphisms $f^i: X^i\ra Y^i$ of $A$-modules such
that $f^id^i_Y=d_X^if^{i+1}$. The morphism $f^\bullet$ is said to be
{\em null homotopic} if there exists a homomorphism $h^i: X^i\ra
Y^{i-1}$ such that $f^i=d_X^ih^{i+1}+h^id_Y^{i-1}$ for all
$i\in\mathbb{Z}$. A complex $X^\bullet$ is called {\em bounded
below} if $X^i=0$ for all but finitely many $i<0$, {\em bounded
above} if $X^i=0$ for all but finitely many $i>0$, and {\em bounded}
if $X^\bullet$ is bounded below and above. We denote by
$\mathscr{C}(A)$ (resp., $\mathscr{C}(\Modcat{A})$) the category of
complexes of finitely generated (resp., all) $A$-modules. The
homotopic category $\mathscr{K}(A)$ is quotient category of
$\mathscr{C}(A)$ modulo the ideals generated by null-homotopic
morphisms. We denote the derived category of $\modcat{A}$ by
$\mathscr{D}(A)$ which is the quotient category of $\mathscr{K}(A)$
with respect to the subcategory of $\mathscr{K}(A)$ consisting of
all the acyclic complexes. The full subcategory of $\mathscr{K}(A)$
and $\mathscr{D}(A)$ consisting of bounded complexes over
$\modcat{A}$ is denoted by ${\mathscr{K}}^b(A)$ and
${\mathscr{D}}^b(A)$, respectively. We denoted by
${\mathscr{C}}^+(A)$ the category of complexes of bounded below, and
by $\mathscr{K}^+(A)$ the homotopic category of $\mathscr{C}^+(A)$.
The full subcategory of $\mathscr{D}(A)$ consisting of bounded below
complexes is denoted by $\mathscr{D}^+(A)$. Similarly, we have the
category $\mathscr{C}^-(A)$ of complexes bounded above, the
homotopic category $\mathscr{K}^-(A)$ of $\mathscr{C}^-(A)$ and the
derived category $\mathscr{D}^-(A)$ of $\mathscr{C}^-(A)$. If we
focus on the category of left $A$-modules, then we have the
homotopic category $\mathscr{K}(\Modcat{A})$ of
$\mathscr{C}(\Modcat{A})$ and the derived category
$\mathscr{D}(\Modcat{A})$ of $\mathscr{C}(\Modcat{A})$.

The two rings $A$ and $B$ are called {\em derived equivalent} if
$\D^b(A)$ and $\D^b(B)$ are equivalent as triangulated categories.
It is equivalent to say that $B$ is isomorphic to
$\End_{\D^b(A)}(T^\bullet)$, where $T^\bullet$ is a complex in
$\K^b(\pmodcat{A})$ satisfying: $(a)$ $T^\bullet$ is
self-orthogonal, that is,
$\Hom_{\K^b(\pmodcat{A})}(T^\bullet,T^\bullet[i])=0$ for all $i\neq
0$, $(b)$ $\add(T^\bullet)$ generates $\K^b(\pmodcat{A})$ as a
triangulated category.


%
%
%
%
%
%
%

In \cite{HX1}, Hu and Xi define the $\mathcal {D}$-split sequences
which occurs in many situation, for instance, the Auslander-Reiten
sequence.

\begin{Def}\label{hw1}\cite{HX1} Let $\mathcal {C}$ be an addictive category and $\mathcal {D}$
a full subcategory of $\mathcal {C}$. A sequence
$$ X \stackrel{f}\lra M \stackrel{g} \lra Y$$
in $\mathcal {C}$ is called a $\mathcal {D}$-split sequence if

$(1)$ $M \in \mathcal {D}$;

$(2)$ f is a left $\mathcal {D}$-approximation of X, and g is a
right $\mathcal {D}$-approximation of Y;

$(3)$ f is a kernel of g, and g is a cokernel of f.
\end{Def}

$\mathcal {D}$-split sequences implies the derived equivalence
between the endomorphism algebras. The following theorem reveals how
to construct derived equivalence from $\mathcal {D}$-split
sequences.

\begin{Lem}\cite{HX1}\label{almost split squence}
Let $\mathcal {C}$ be an additive category and $M$ an object in
$\mathcal {C}$. Suppose

$$X\stackrel{f}\ra M'\stackrel{g}\ra Y$$
is an almost $\add (M)$-split sequence in $\mathcal {C}$. Then the
endomorphism ring $\End_\mathcal {C}(X\oplus M)$ of $X\oplus M$ and
the endomorphism ring $\End_\mathcal {C}(Y\oplus M)$ of $Y\oplus M$
are derived equivalent.
\end{Lem}

\section{Results and proofs}
To prove our results, we first establish a fact.

\begin{Lem}\label{144}  Let R be a ring with identity.

{\rm (1)}Let M be a neotherian left R-module, and let $f:M \ra M$ be
a surjective homomorphism, then f is injective.

{\rm (2)}Let M be an artinian left R-module, and let $f: M \ra M$ be
an injective homomorphism, then f is surjective.

 \label{lem1}
\end{Lem}

{\it Proof.} We only prove the first part of the lemma. The second
part of the lemma is similar. Set $f^{k}$=$ \overbrace{f\cdots
f}\limits^k$. Note that $\Ker f^k$ are submodules of $M$ and $\Ker
f^i$ is the submodule of $\Ker f^{i+1}$ for any $i\geq 1$. Since $M$
is a neotherian module, there exists $i_0\geq 1$, satisfying $\Ker
f^{i_0}=\Ker f^{i_0+1}$. Then we have the following commutative
diagram:
 $$\xymatrix{
  0 \ar[r] & \Ker f^{i_0}\ar[d]_{id} \ar[r] & M \ar[d]_{id}
  \ar[r]^{^{f^{i_0}}} & M \ar[d]_{^{f^{i_0}}} \ar[r]& 0 \\
  0 \ar[r] & \Ker f^{i_0+1} \ar[r] & M\ar[r]^{f^{i_0+1}} & M\ar[r] &
  0}$$ By the snake lemma, we can get $\Ker f=0$. So f is injective.
  $\square$

Now, let us prove the main result in this paper.

{\bf Proof of Theorem \ref{Theo}:} Set $\Gamma=M_n(A)$, the $n
\times n$ matrix over $A$.
Denote by $e_i$ the matrix which has $1_A$ in the $(i,i)$-th
position and zeros elsewhere for $1\leq i\leq n$. So $e_1, e_2,
\cdots, e_n$ are piecewise orthogonal idempotents in $\Lambda$, such
that $1_{\Lambda} = e_1 + e_2 + \cdots + e_n$.

Since $\Lambda$ is a subring of $\Gamma$ with the same identity, the
ring $\Gamma$ can be considered as a $\Lambda$-module just by
restriction of the scalars of $\Gamma$ to $\Lambda$.

Now, we consider the exact sequence
$$0\lra \Lambda\stackrel{\lambda}\lra \Gamma \stackrel{\pi}\lra L\lra 0 $$
in $\Lambda$-Mod, where $\lambda$ is the inclusion map and $L$ is
the cokernel of $\lambda$. To show the Theorem, we prove the
following statements.

$(1)$ The sequence
$$0\lra \Lambda \stackrel{\lambda}\lra \Gamma \stackrel{\pi}
\lra L \lra 0$$ is an almost $\add(\Lambda e_1)$- split sequence in
$\Lambda$-Mod.

In fact, we shall check that all conditions in Definition \ref{hw1}
are satisfied.

Since the left $\Lambda$-module $_{\Lambda}\Gamma$ is a direct sum
of some copies of $\Lambda e_1$, we have $\Gamma \in$ $\add(\Lambda
e_1)$. Clearly, $\Lambda e_1$ is projective as a left
$\Lambda$-module, then we have an exact sequence
$$0\lra \Hom_{\Lambda}(D, \Lambda) \stackrel{(-, \lambda)} \lra
\Hom_{\Lambda}(D, \Gamma) \stackrel{(-, \pi)} \lra \Hom_{\Lambda}(D,
L) \lra 0$$ for any $D \in$ $\add(\Lambda e_1)$.

This means that the homomorphism $\pi : {}_{\Lambda}\Gamma\lra{}
_{\Lambda}\Gamma$ is a right $\add(\Lambda e_1)$-approximation of
$_{\Lambda}L$. Now, we prove that the homomorphism $\lambda :{}
_{\Lambda}\Lambda \lra {}_{\Lambda} \Gamma$ is a left $\add(\Lambda
e_1)$-approximation of $\Lambda$. In fact, every left
${\Lambda}$-module homomorphism $g :{}\Lambda \lra\Lambda e_1$ is
determined by $g(1)$, the image of 1 under $g$. Similarly, every
left $\Gamma$-module homomorphism $h :{} \Gamma \ra\Gamma e_1$ is
determined by $h(1)$, the image of 1 under $h$. Note that $\Gamma
e_1$ and $\Lambda e_1$ are isomorphic as left $\Lambda$-module, and
any left $\Gamma$-module homomorphism is also left $\Lambda$-module
homomorphism. So we assume that $g_1 :{} \Gamma \ra \Lambda e_1$ is
a left $\Lambda$-module homomorphism which sends 1 to $g(1)$. Then
the homomorphism $g_1$ satisfies $g = i g_1$. Thus we have proved
that the homomorphism $\lambda :{} _\Lambda \Lambda \ra{}_\Lambda
\Gamma$ is a left $\add(\Lambda e_1)$-approximation. Hence $(1)$ is
proved.

Note that the sequence
$$0\lra \Lambda e_i \stackrel{\lambda_i}\lra \Lambda e_1 \stackrel{\pi_i}
\lra L_i \lra 0  \quad\quad (*)  $$ where $\lambda_i$ is inclusion
map and $L_i$ is the cokernel of $\lambda_i$, are almost
$\add(\Lambda e_1)$-split sequences for $2\leq i\leq n$.

By Lemma \ref{almost split squence}, the ring $\Lambda$ and the
endomorphism ring $\End_{\Lambda}(L_2 \oplus L_3 \oplus \cdots
\oplus L_n \oplus \Lambda e_1)$ are derived equivalent via a tilting
module $L_2 \oplus L_3 \oplus \cdots \oplus L_n \oplus \Lambda e_1$.

$(2)$ The ring $\Sigma$ and the endomorphism ring
$\End_{\Lambda}(L_2 \oplus L_3 \oplus \cdots \oplus L_n \oplus
\Lambda e_1)$ are isomorphic as rings.

Indeed, we note that
$$\End_{\Lambda}(L_2 \oplus L_3 \oplus \cdots \oplus L_n \oplus
\Lambda e_1) \cong\begin{pmatrix}
(L_2,L_2) & (L_2, L_3) & \cdots & (L_2, \Lambda e_1)\\
(L_3, L_2 ) & (L_3,L_3) & & (L_3, \Lambda e_1)\\
\vdots & & \ddots & \vdots \\
(\Lambda e_1, L_2) & (\Lambda e_1, L_3) & \cdots & (\Lambda
e_1,\Lambda e_1)
\end{pmatrix}$$
as rings.\medskip

In the following, we calculate the endomorphism ring
$\End_{\Lambda}(L_2 \oplus L_3 \oplus \cdots \oplus L_n \oplus
\Lambda e_1)$.

The morphism set $\Hom_\Lambda(L_i, \Lambda e_i)=0$ for $2\leq i\leq
n$. Applying the functor $\Hom_\Lambda(-,\Lambda e_1)$ to the exact
sequence $0\ra \Lambda e_i\ra \Lambda e_1\ra L_i\ra 0$, we can get
the following exact sequence
$$0\ra \Hom_\Lambda(L_i, \Lambda e_1)\ra \Hom_\Lambda(\Lambda e_1, \Lambda e_1)\ra \Hom_\Lambda(\Lambda e_i, \Lambda e_1)\ra 0$$
in $\mathbb{Z}$-Mod for $2\leq i\leq n$.

Note that both of $\Hom_\Lambda(\Lambda e_1,\Lambda e_1)$ and
$\Hom_\Lambda(\Lambda e_i, \Lambda e_1)$ are isomorphic to $A$ in
$\mathbb{Z}$-Mod. By Lemma \ref{144}, we have
$\Hom_\Lambda(L_i,\Lambda e_1)=0$ for $2\leq i\leq n$.

\medskip

For simplicity, we denote the set $e_i\Lambda e_j$ by $\Lambda_{ij}$
for $1\leq i,j\leq n$. The morphism $\mu_x$ denote the right
multiplication by $x$.

Let $b$ be an element in $\Lambda_{ij}$ for $2\leq i,j\leq n$. Since
$\lambda_i$ is the left $\Lambda e_1$-approximation of $\Lambda
e_i$, we have a morphism $\mu_a: \Lambda e_1\ra \Lambda e_1$ such
that $\lambda_i\mu_a=\mu_b\lambda_j$ where $a$ is an element in $A$.
Thus we can get an element $\alpha_b$ in $\Hom_\Lambda(L_i,L_j)$
such that $\pi_i\alpha_b=\mu_a\pi_j$. It follows from the
commutativity of the left square that $a=b$. For a given morphism
$\mu_b$, there is a unique $\alpha_b$ satisfying
$\pi_i\alpha_b=\mu_b\pi_j$. Note that $\pi_i\mu_b=\mu_b\pi_j$, we
can get $\alpha_b=\mu_b$.
$$\xymatrix{
  0  \ar[r] & \Lambda e_i \ar[d]_{\mu_b} \ar[r]^{\lambda_i} & \Lambda e_1 \ar@{-->}[d]_{\mu_a} \ar[r]^{\pi_i} & L_i \ar@{-->}[d]_{\alpha_b} \ar[r]^{} & 0 \\
 0 \ar[r] & \Lambda e_j \ar[r]^{\lambda_j} & \Lambda e_1 \ar[r]^{\pi_j} & L_j \ar[r]^{} & 0   }\quad\quad (*)$$
where $\lambda_i, \lambda_j$ are the inclusion maps and $L_i, L_j$
are the cokernels of $\lambda_i$ and $\lambda_j$, respectively.

Thus, we can define a set of maps from $\Lambda_{ij}$ to
$\Hom_\Lambda(L_i,L_j)$ for $2\leq i,j\leq n$.

Define
$$\phi_{ij}: \Lambda_{ij}\ra \Hom_\Lambda(L_i, L_j)$$
$$b\mapsto \alpha_b.$$
for $2\leq i,j\leq n$.

$(a)$ The map $\phi_{ij}$ is well-defined.
$$\xymatrix{
  0  \ar[r] & \Lambda e_i \ar[d]_{\mu_b} \ar[r]^{\lambda_i} & \Lambda e_1 \ar[d]_{\mu_b} \ar[r]^{\pi_i} & L_i \ar@{-->}[ld]_{s_i}\ar[d]_{\alpha_b} \ar[r]^{} & 0 \\
 0 \ar[r] & \Lambda e_j \ar[r]^{\lambda_j} & \Lambda e_1 \ar[r]^{\pi_j} & L_j \ar[r]^{} & 0   }$$

Suppose that $b=0$, we have $\lambda_i\mu_b=0$. It follows that
there is a morphism $s_i: L_i\ra \Lambda e_1$ such that $\mu_b=\pi_i
s_i$. Thus we have $\pi_i\alpha_b=\mu_b\pi_j=\pi_is_i\pi_j$. Since
$\pi$ is surjective, we have $\alpha_b=s_i\pi_j$. Note that
$\Hom_\Lambda(L_i, \Lambda e_1)=0$ for $2\leq i\leq n$, we obtain
$\alpha_b=0$. Hence $\phi_{ij}$ is well-defined.

\medskip

$(b)$ The morphism $\phi_{ij}$ is surjective.

Let $\alpha$ be an element in $\Hom_\Lambda(L_i,L_j)$. Note that
$\Lambda e_1$ is projective module over $\Lambda$, thus there exists
a morphism $\mu_{a}: \Lambda e_1\ra \Lambda e_1$ such that
$\mu_a\pi_j=\pi_i\alpha$ where $a$ is an element in $\Lambda_{11}$.
Thus there is a unique morphism $\mu_{b}: \Lambda e_i\ra \Lambda
e_j$ such that $\lambda_i\mu_b=\mu_b\lambda_j$ for $b\in
\Lambda_{ij}$. So $\phi_{ij}$ is surjective.

\medskip

$(c)$ The description of $\Ker\phi_{ij}$, i.e., $\Ker\phi_{ij}=I_j$
for $2\leq i, j\leq n$.
$$\xymatrix{
  0  \ar[r] & \Lambda e_i \ar[d]_{\mu_b} \ar[r]^{\lambda_i} & \Lambda e_1 \ar@{-->}[ld]_{t_i}\ar[d]_{\mu_b} \ar[r]^{\pi_i} & L_i \ar[d]_{\alpha_b} \ar[r]^{} & 0 \\
 0 \ar[r] & \Lambda e_j \ar[r]^{\lambda_j} & \Lambda e_1 \ar[r]^{\pi_j} & L_j \ar[r]^{} & 0.   }$$

Suppose that $\alpha_b=0$. Then we have $\mu_b\pi_j=0$. So there is
a morphism $t_i: \Lambda e_1\ra\Lambda e_j$ such that
$\mu_b=t_i\lambda_j$. Note that there exist $c\in I_j$, $d\in A$
such that $t_i=\mu_c$, $\lambda_j=\mu_d$. Thus, we have $b=cd$,
i.e., $b\in I_j$. Conversely, suppose that $b$ is an element in
$I_j$, we have $\mu_b\pi_j=\pi_i\alpha_b=0$. It follows that
$\alpha_b=0$. Hence, $\Ker \phi_{ij}=I_j$ for $2\leq i, j\leq n$.

\medskip

$(d)$ $\phi_{ij}$ preserves addition and multiplication.

It is easy to prove that $\phi_{ij}$ preserves addition.

Now, we turn to prove that $\phi_{ij}$ preserves multiplication,
i.e., $\phi_{ij}(b)\phi_{jk}(b')=\phi_{ik}(bb')$ for $2\leq i, j\leq
n$ where $b$ and $b'$ are elements of $\Lambda_{ij}$ and
$\Lambda_{jk}$ respectively.

It suffices to prove that $\alpha_{bb'}=\alpha_b\alpha_{b'}$. Since
$\mu_b\mu_{b'}=\mu_{bb'}$, we have
$\lambda_i\mu_{aa'}=\lambda_i\mu_c$ where $\mu_c$ is a morphism
induced by $\mu_{bb'}$. Thus $\mu_{aa'}-\mu_c$ factorizes through
$\pi_i$. Note that $\Hom_{\Lambda}(L_i,\Lambda e_1)=0$ for $2\leq
i\leq n$, we have $\mu_{aa'}=\mu_c$. Hence we get
$\alpha_{bb'}=\alpha_b\alpha_{b'}$.

\medskip

Now, we can define a map
$$\phi=(\phi_{ij}):\quad
(e_2+\cdots+e_n)\Lambda(e_2+\cdots+e_n)\ra \small\begin{pmatrix}
(L_2, L_2) & (L_2, L_3) & \cdots &(L_2, L_n)\\
(L_3, L_2 ) & (L_3, L3) & &(L_3, L_n)\\
\vdots &\vdots & \ddots&\vdots\\
(L_n, L_2)&(L_n, L_3)&\cdots&(L_n, L_n)
\end{pmatrix}$$

$$(a_{ij})_{i-1,j-1}\mapsto (\phi_{ij}(a_{ij}))_{i-1,j-1}$$
for $2\leq i, j\leq n$.

The map $\phi_{ij}$ is well-defined and surjective, so is the map
$\phi$. It follows from that $\phi_{ij}$ preserves addition and
multiplication for $2\leq i, j\leq n$ that $\phi$ is a ring
homomorphism. The kernel of $\phi$ is
$$\small\begin{pmatrix}
I_2&I_3&\cdots&I_n\\
I_2&I_3&\cdots&I_n\\
\vdots&\vdots&&\vdots\\
I_2&I_3&\cdots&I_n
\end{pmatrix}.$$

Thus, we have a ring isomorphism
$$\overline{\phi}: \small \begin{pmatrix}
A_2/I_2&0&\cdots&&0\\
I_{32}/I_2&A_3/I_3&\ddots&&\vdots\\
I_{42}/I_2&I_{43}/I_{3}&A_4/I_4&&\\
\vdots&\vdots&\vdots&\ddots&0\\
I_{n2}/I_2&I_{n3}/I_3&I_{n4}/I_4&\cdots&A_n/I_n
\end{pmatrix}\lra\small\begin{pmatrix}
(L_2, L_2) & (L_2, L_3) & \cdots &(L_2, L_n)\\
(L_3, L_2 ) & (L_3, L_3) & &(L_3, L_n)\\
\vdots &\vdots & \ddots&\vdots\\
(L_n, L_2)&(L_n, L_3)&\cdots&(L_n, L_n)
\end{pmatrix}.$$

Meanwhile, we have algebra isomorphism $\varphi_2$
$$\varphi_2: A\ra \End_\Lambda(\Lambda e_1)$$
$$a\mapsto (f_a:\lambda e_1\mapsto\lambda ae_1)$$
and isomorphism of abelian groups $\varphi_3$
$$\varphi_3: (A/I_2, A/I_3, \cdots,A/I_n)\ra (\Hom_\Lambda(\Lambda e_1,L_2),\Hom_\Lambda(\Lambda e_1,L_3),\cdots,\Hom_\Lambda(\Lambda e_1,L_n))$$
$$(m_1,m_2,\cdots,m_n)\mapsto (f_{m_2},f_{m_3},\cdots,f_{m_n})$$
where $f_{m_i}: \lambda e_1\mapsto \lambda m_ie_1$ for $2\leq i\leq
n$.

Now, Set
$$\varphi=\begin{pmatrix}
  \overline{\phi}&0\\
  \varphi_3&\varphi_2
\end{pmatrix}: \Sigma\ra \begin{pmatrix}
\End_{\Lambda}(L_2) & \Hom_{\Lambda}(L_2, L_3) & \cdots & \Hom_{\Lambda}(L_2, \Lambda e_1)\\
\Hom_{\Lambda}(L_3, L_2 ) & \End_{\Lambda}(L_3) & & \Hom_{\Lambda}(L_3, \Lambda e_1)\\
\vdots & & \ddots & \vdots \\
\Hom_{\Lambda}(\Lambda e_1, L_2) & \Hom_{\Lambda}(\Lambda e_1, L_3)
& \cdots & \End_{\Lambda}(\Lambda e_1)
\end{pmatrix}$$
$$\small\begin{pmatrix}
  r_{22}&r_{23}&\cdots&r_{2n}&0\\
  r_{32}&r_{33}&\cdots&r_{3n}&\vdots\\
  \vdots&\vdots&\ddots&\vdots&\\
  r_{n2}&r_{n3}&\cdots&r_{nn}&0\\
  m_2&m_3&\cdots&m_n&a
\end{pmatrix}\mapsto\small\begin{pmatrix}
  \overline{\phi_{22}}(r_{22})&\overline{\phi_{23}}(r_{23})&\cdots&\overline{\phi_{2n}}(r_{2n})&0\\
  \overline{\phi_{32}}(r_{32})&\overline{\phi_{33}}(r_{33})&\cdots&\overline{\phi_{3n}}(r_{3n})&\vdots\\
  \vdots&\vdots&\ddots&\vdots&\\
  \overline{\phi_{n2}}(r_{n2})&\overline{\phi_{n3}}(r_{n3})&\cdots&\overline{\phi_{nn}}(r_{nn})&0\\
  \varphi_3(m_2)&\varphi_3(m_3)&\cdots&\varphi_3(m_n)&\varphi_2(a)
\end{pmatrix}.$$

Clearly, the map $\varphi$ is an isomorphism of abelian groups. And
it easy to  check that $\varphi$ is a ring isomorphism. The proof is
completed. $\square$



\medskip

As a direct consequence of Theorem \ref{Theo}, we have the following
corollary.

\begin{Koro}Let $A$ be a ring with identity, and suppose that $I_2, I_3, \cdots, I_n$ are
ideals of $A$.

$(1)$ The two rings
$$\footnotesize\begin{pmatrix}
  A & I_2& I_3& \cdots & \cdots & I_n\\
  A & A & I_3 & \cdots & \cdots & I_n\\
  A & I_2 & A &I_4 & \cdots & I_n\\
  A & I_2 & I_3 & A & \cdots & I_n\\
  \vdots & \vdots & \vdots && \ddots & \vdots \\
  A &I_2 & I_3 & \cdots & I_{n-1} & I_n
\end{pmatrix}\mbox{ and }\footnotesize\begin{pmatrix}
  A/I_2 &0 & \cdots& & & 0\\
  0& A/I_3 &\ddots & & &\vdots\\
  \vdots& \ddots& A/I_4 & & &\\
  & & & \ddots & &\\
  0&\cdots & &0 & A/I_n &0\\
  A/I_2 & A/I_3 & A/I_4 & \cdots & A/I_n & A
\end{pmatrix} $$
are derived equivalent.

$(2)$ The two rings
$$\footnotesize\begin{pmatrix}
  A & I_2 & I_3 &  \cdots & I_{n-1} & I_n\\
  A & A & I_3 &  \cdots & I_{n-1} & I_n\\
  A & A & A &  \cdots & I_{n-1} & I_n\\
  \vdots & \vdots & \vdots &  \ddots & I_{n-1} & I_n\\
  \vdots & \vdots & \vdots &   & A & I_n\\
  A & A & A & A & \ldots & A
\end{pmatrix}\mbox{ and }\footnotesize\begin{pmatrix}
  A/I_2 & 0& \cdots& & &  0  \\
  A/I_2 & A/I_3 &\ddots & &  \\
  A/I_2 & A/I_3 & A/I_4 & & &\vdots\\
  \vdots & \vdots & \vdots &\ddots& &\\
  \vdots & \vdots & \vdots && A/I_n &0\\
  A/I_2 & A/I_3 & A/I_4 & \cdots & A/I_n &A
\end{pmatrix}$$
are derived equivalent. \label{cor}
\end{Koro}

\section{Applications}

In section 3, we have constructed derived equivalences of matrix
subrings. In this section, we will give some applications of the
main result. At first, we define a class of rings which are called
general block extensions of rings with respect to the decomposition of identity.
Then we calculate the finitistic dimension of a general block
extension of ring with respect to the decomposition of identity. At last, we
will consider the finiteness of finitistic dimension of a tiled
triangular ring. Proposition \ref{fd}, will give a condition under
which the finitistic dimensions of tiled triangular algebras are
finite.

The following lemmas, which are taken from \cite{PX, K}, are useful
for this section.

\begin{Lem}\cite{PX}\label{PX} If two left coherent rings $A$ and $B$ are derived equivalent, and if $T^{\bullet}$ is a
tilting complex over $A$ with $n+1$ non-zero terms such that $B\cong
\End(T^{\bullet})$, then $\fd (A)-n\leq \fd (B)\leq \fd (A)+n$.

\end{Lem}

\noindent{\bf Remark.}\label{px1} In Lemma \ref{PX}, if we replace
$A$ and $B$ by arbitrary ring with identity, then $\lfd (A)-n\leq
\lfd (B)\leq \lfd (A)+n$. The proof is similar.

\begin{Lem}\cite{K}\label{comm} Let $A$ be an ring with
identity, and let $T^{\bullet}$ be a tilting complex over $A$ with
$\End(T^{\bullet})\simeq B$. If $T^{\bullet}$ has $n+1$ non-zero
terms, where $n\geq 0$, then the following statements hold.

$(a)$ $\lgd(A)-n\leq \lgd(B)\leq \lgd(A)+n$;

$(b)$ $\injdim(_AA)-n\leq \injdim(_BB)\leq \injdim(_AA)+n
$.\label{lem4}
\end{Lem}

The following lemma about the estimation of global dimension and
finitistic dimension can be found in \cite[Corollary 4.21,
p.70]{FGR}.

\begin{Lem}\cite{FGR}\label{456} Let R and S be rings. Let M be an S-R bimodule and $\Lambda:=
\begin{pmatrix}
  R & 0\\
  M & S
\end{pmatrix}$.
Then the following inequalities hold:

$(1)$ $\lfd(S)\leq \lfd(\Lambda)\leq 1+\lfd(R) +\lfd(S)$.

$(2)$ $\lfd(\Lambda)\geq sup\{ pd(_RA)\leq \infty \mid
A\in\Modcat{R} \mbox{ satisfying }Tor^A_i(M,A) = 0 \mbox{ for all }
i \} $. If $M$ is flat as a right-$R$-module, then
$\lfd(\Lambda)\geq \lfd(R)$.

$(3)$ If $\pd(_SM)\leq \infty$, then $\pd(_SM)+1\leq \lfd(\Lambda)$
$ \leq max\{\lfd(R)+\pd(_SM)+1, \lfd(S)\}$.

$(4)$ $max\{\lgd(R), lgl.dim(S), \pd(_SM)+1\}\leq \lgd(\Lambda)\leq
max\{\lgd (R)+\pd(_SM)+1, \lgd(S)\}$.
 \label{lem3}
\end{Lem}

The corresponding statements hold for the right homological
dimensions over $\Lambda$.



\subsection{General block extensions of rings with respect to the decomposition of identity}

In this part, we will define a class of rings which contains
hereditary orders, block extensions of basic algebras.

\begin{Def}\label{block extension}
Let $A$ be a ring with identity $1_A$. And $1_A=e_1+e_2+\cdots+e_m$
is a decomposition of the identity where $e_i$ is idempotent. Then
$A$ can be represented as the following matrix form
$$A=\small\begin{pmatrix}
  e_1Ae_1&e_1Ae_2&\cdots&e_1Ae_m\\
  e_2Ae_1&e_2Ae_2&\cdots&e_2Ae_m\\
  \vdots&\vdots&\ddots&\vdots\\
  e_mAe_1&e_mAe_2&\cdots&e_mAe_m
\end{pmatrix}.$$
Set $A_i=e_iAe_i$ and $A_{ij}=e_iAe_j$. Then $A_i$ is the subring of
$A$ with identity element $e_i$, and $A_{ij}$ is a $(A_i,
A_j)$-bimodule.

Let $n_1, n_2, \cdots, n_m\in \mathbb{N}$. For $1\leq i, s\leq m,
1\leq j\leq n_i$ and $1\leq t\leq n_s$, we define
$$P=A(n_1, n_2, \cdots, n_m)=\small\begin{pmatrix}
  P(1,1)&P(1,2)&\cdots&P(1,m)\\
  P(2,1)&P(2,2)&\cdots&P(2,m)\\
  \vdots&\vdots&\ddots&\vdots\\
  P(m,1)&P(m,2)&\cdots&P(m,m)
\end{pmatrix}$$
which is contained in the ring
$\End_A((Ae_1)^{n_1}\oplus\cdots\oplus(Ae_m)^{n_m})$ with the
restrictions of the binary operations of addition and multiplication
of $\End_A((Ae_1)^{n_1}\oplus\cdots\oplus(Ae_m)^{n_m})$.
$$P(i,s)=\small\begin{pmatrix}P_{i1,s1}&P_{i1,s2}&\cdots&P_{i1,sn_s}\\
P_{i2,s1}&P_{i2,s2}&\cdots&P_{i2,sn_s}\\
\vdots&\vdots&\ddots&\vdots\\
P_{in_i,s1}&P_{in_i,s2}\cdots&\cdots&P_{in_i,sn_s}
\end{pmatrix}_{n_i\times n_s}$$
satisfies that $P_{ip,sq}$ is a $(A_i,A_s)$-bimodule.

For $P(i,s)$, there are three cases:

Case $I: i=s$.
$$P(i,s):=\footnotesize\begin{pmatrix}
A_i & I_{i2} &I_{i3} & \cdots & I_{i(n-1)} & I_{in}\\
A_i & B_{i2} & I_{i3} & \cdots & I_{i(n-1)} & I_{in}\\
A_i & I_{i32}& B_{i3}& \cdots & I_{i(n-1)} & I_{in}\\
A_i & I_{i42} & I_{i43}& B_{i4} & \cdots & I_{in}\\
\vdots & \vdots &\ddots &\ddots  & \ddots& \vdots\\
A_i & I_{in_i2} & \cdots & \cdots & I_{in_i(n_i-1)} & B_{in_i}
 \end{pmatrix}$$
 where $B_{il}$ is the subring of $A_i$ with the same
 identity, $I_{il}$ and $I_{ipq}$ are ideals of $A_i$ satisfying
 $I_{in_i}\subseteq I_{i(n_i-1)}\subseteq\cdots\subseteq I_{i2}, I_{il}\subseteq I_{ipq}, I_{il}\subseteq B_{il}$
 for $2\leq l\leq n_i, 3\leq p\leq n_i, 2\leq q\leq n_i-1$.

Case $II: i<s$.
$$P(i,s):=\small\begin{pmatrix}
 A_{is}&P_{i1,s2}&\cdots&P_{i1,sn_s}\\
 A_{is}&P_{i1,s2}&\cdots&P_{i1,sn_s}\\
\vdots&\vdots&\vdots&\vdots\\
A_{is}&P_{i1,s2}&\cdots&P_{i1,sn_s}
\end{pmatrix}$$

Case $III: i>s$.
$$P(i,s):=\small\begin{pmatrix}
  A_{is}&P_{i1,s2}&\cdots&P_{i1,sn_s}\\
  A_{is}&P_{i2,s2}&\cdots&P_{i2,sn_s}\\
  \vdots&\vdots&\ddots&\vdots\\
  A_{is}&P_{in_i,s2}&\cdots&P_{in_i,sn_s}
\end{pmatrix}.$$

Suppose that $\sum^{m}_{l=1}P(i, l)P(l, j)\subseteq P(i, j)$. Then
$P$ is a ring and called  general block extension of $A$ with respect to the
decomposition of identity.
\end{Def}

General block extension of $A$ with respect to the decomposition of identity
contains many classes of subrings of $M_n(A)$. In the following, we
will give some examples.

\noindent{\bf Example.} $(1)$ In Definition \ref{block extension},
we assume that $m=1$. Then
$$P=\small\begin{pmatrix}
A & I_2 &I_3 & \cdots & I_{n-1} & I_n\\
A & A_2 & I_3 & \cdots & I_{n-1} & I_n\\
A & I_{32}& A_3& \cdots & I_{n-1} & I_n\\
A & I_{42} & I_{43}& A_4 & \cdots & I_n\\
\vdots & \vdots &\ddots &\ddots  & \ddots& \vdots\\
A & I_{n2} & \cdots & \cdots & I_{n,n-1} & A_n
 \end{pmatrix}$$
where $A$ is a ring with identity,
 $A_i$ are subrings of $A$ with the same identity
for $2\leq i \leq n$. $I_i, I_{ij}$ are ideals of $A$ for $2\leq
i\leq n, 2\leq j\leq n-1$. In particular, set $I_i=a\Omega,\,
A=A_i=\Omega$ for $2\leq i\leq n$ and $I_{ij}=\Omega$ for $2\leq
j<i\leq n$, where $\Omega$ is a local $R$-order and $a$ is a regular
element in $\Omega$. Then $\Omega/a\cdot\Omega$ is local and $P$ is
a QH-order with associated ideal $J=\omega\cdot P$, where
$$\omega=\small\begin{pmatrix}
  0&1&0&\cdots&0\\
  0&0&1&\cdots&0\\
  &&\cdots&&\\
  0&0&\cdots&1&0\\
  0&0&0&\cdots&1\\
  a&0&\cdots&\cdots&0
\end{pmatrix}_{n\times n}.$$


$(2)$ In Definition \ref{block extension}, let $A$ be a basic
algebra, and let $\{e_1, \cdots, e_m\}$ be a complete set of
orthogonal primitive idempotents of $A$.

set
$$P(i,s)=\small\begin{pmatrix}P_{i1,s1}&P_{i1,s2}&\cdots&P_{i1,sn_s}\\
P_{i2,s1}&P_{i2,s2}&\cdots&P_{i2,sn_s}\\
\vdots&\vdots&\ddots&\vdots\\
P_{in_i,s1}&P_{in_i,s2}\cdots&P_{in_i,sn_s}&P_{in_i, sn_s}
\end{pmatrix}=
\left\{\begin{array}{ll}
  \begin{pmatrix}
    A_i&\cdots&A_i\\
    &\ddots&\vdots\\
    \rad(A_i)&&A_i
  \end{pmatrix}&(i=s)\\
  \begin{pmatrix}
    A_{is}&\cdots&A_{is}\\
    \vdots&&\vdots\\
    A_{is}&\cdots&A_{is}
  \end{pmatrix}&(i\neq s)
\end{array}\right.$$
Then $P$ is called the {\em block extension} of $A$ which can be
found in \cite{O3}. In particular, if $A$ is a basic QF-algebra,
then $P$ is a basic Harada algebra (see \cite{O3}).

%
%

\begin{Theo}\label{theo5}Let $A$ be a ring with identity.
Let $P$ be a general block extension of $A$ with respect to the decomposition of
identity.
Then
$$\lfd (A)-1\leq \lfd (P)\leq \lfd (A)+\sum^{m}_{j=1}\sum^{n_j}_{i=2}\lfd (A_j/I_{ji})+\sum^m_{i=1}n_i-m.$$
\end{Theo}

{\bf Proof.} Denote by $e_{\sum^{i}_{l=1}n_l+j}$ the matrix which
has $1_{A_{i+1}}$ in the $(\sum^{i}_{l=1}n_l+j$,
$\sum^{i}_{l=1}n_l+j)$-th\linebreak position and zeros elsewhere for
$1\leq j\leq n_{i+1}$, $0\leq i\leq m-1$. Thus $e_1,\cdots, e_{n_1},
e_{n_1+1},\cdots, e_{n_1+n_2},\linebreak \cdots,
e_{\sum^{m-1}_{l=1}n_l+1}, \cdots, e_{\sum^m_{l=1}n_l}$ are
piecewise orthogonal idempotents in $P$ such that $1_P = e_1 +
\cdots + e_{\sum^m_{l=1}n_l}$.

Set $\Upsilon=\End_P(\overbrace{Pe_1\oplus\cdots\oplus
Pe_1}\limits^{n_1}\oplus\cdots\oplus\overbrace{Pe_{\sum^{m-1}_{l=1}n_l1}\oplus\cdots\oplus
Pe_{\sum^{m-1}_{l=1}n_l1}}\limits^{n_m})$.
Since $P$ is a subring of $\Upsilon$ with the same identity,
$\Upsilon$ can be viewed as a $P$-module by restriction of the
scalars of $\Upsilon$ to $P$. There is an exact sequence
$$0\ra P\stackrel{\lambda}\ra \Upsilon\stackrel{\pi}\ra L\ra 0$$
in $\Modcat{\Upsilon}$, where $\lambda$ is the inclusion map and $L$
is the cokernel of $\lambda$. By Theorem \ref{almost split squence},
two rings $P$ and $\End_P((\oplus_{i=1,\cdots,m\atop
j=2,\cdots,m}L_{n_in_j})\oplus(\oplus^{\sum^m_{i=1}n_i}_{i=1}Pe_i))$
are derived equivalent.

Note that $\End_P((L_{n_12}\oplus\cdots\oplus
L_{n_1n_1})\oplus\cdots\oplus(L_{n_m1}\oplus\cdots\oplus
L_{n_mn_m})\oplus(Pe_{\sum^{m-1}_{l=1}n_l+1}\oplus\cdots\oplus
Pe_1))\cong$
$$\footnotesize\begin{pmatrix}
((L_{n_1,i}, L_{n_1,j}))_{2\leq i,j\leq n_1}&((L_{n_1i},
L_{n_2j}))_{2\leq i\leq n_1\atop 2\leq j\leq n_2}&\cdots&((L_{n_1i},
L_{n_mj}))_{2\leq i\leq n_1 \atop 2\leq j\leq n_m}&((L_{n_1i},
Pe_i))_{2\leq i\leq n_1\atop 1\leq j\leq m}\\
((L_{n_2,i}, L_{n_1,j}))_{2\leq i\leq n_2\atop 2\leq j\leq
n_1}&((L_{n_2i}, L_{n_2,i}))_{2\leq i,j\leq n_2}&\cdots&((L_{n_2i},
L_{n_mj}))_{2\leq i\leq n_2 \atop 2\leq j\leq n_m}&((L_{n_2i},
Pe_j))_{2\leq i\leq n_1\atop 1\leq j\leq m}\\
\vdots&\vdots&\ddots&\vdots&\vdots\\
((L_{n_m,i}, L_{n_1,j}))_{2\leq i\leq n_m\atop 2\leq j\leq
n_1}&((L_{n_mi}, L_{n_2,j}))_{2\leq i\leq n_m\atop 2\leq j\leq
n_2}&\cdots&((L_{n_mi},L_{n_mj}))_{2\leq i,j\leq n_m}&((L_{n_mi},Pe_j))_{2\leq i\leq n_m\atop 1\leq j\leq m}\\
((Pe_i, L_{n_1j}))_{1\leq i\leq n_m\atop 2\leq j\leq n_1}&((Pe_i,
L_{n_2j}))_{2\leq i\leq n_m\atop 2\leq j\leq
n_2}&\cdots&((L_{n_mi},L_{n_mj}))_{2\leq i,j\leq
n_m}&((L_{n_mi},Pe_j))_{2\leq i\leq n_m\atop 1\leq j\leq m}
\end{pmatrix}$$

$(1)$ $\Hom_P(L_{n_pi}, L_{n_qj})=0$ for $ 2\leq i\leq n_p, 2\leq
j\leq n_q, 1\leq p<q\leq m$.

There are exact sequences
$$0\ra Pe_{\sum^{p-1}_{l=1}n_l+i}\ra Pe_{\sum^{p-1}_{l=1}n_l+1}\ra L_{n_{p}i}\ra 0.\quad\quad(**)$$

Applying the functor $\Hom_P(-, L_{n_{q}j})$ to $(**)$,
we get an exact sequence
$$0\ra \Hom_P(L_{n_pi}, L_{n_qj})\ra \Hom_P(Pe_{\sum^{p-1}_{l=1}n_l+1}, L_{n_qj})\ra \Hom_P(Pe_{\sum^{p-1}_{l=1}n_l+i}, L_{qj})\ra 0.$$

By calculation, we have $\Hom_P(Pe_{\sum^{p-1}_{l=1}n_l+1},
L_{n_qj})=\Hom_P(Pe_{\sum^{p-1}_{l=1}n_l+i}, L_{qj})=0$. Thus,
$\Hom_P(L_{n_{p}i}, L_{n_{q}j})=0$ for $2\leq i\leq n_p, 2\leq j\leq
n_q, 1\leq p<q\leq m.$

$(2)$ $\Hom_P(L_{n_ki}, L_{n_kj})=0$ for $2\leq i<j\leq n_k, 1\leq
k\leq m$.

Apply the functor $\Hom_P(-,L_{n_kj})$ to the exact sequence
$$0\ra Pe_{\sum^{k-1}_{j=1}n_j+i}\ra Pe_{\sum^{k-1}_{j=1}n_j+1}\ra L_{n_{k}i}\ra 0.$$
We have an exact sequence
$$0\ra \Hom_P(L_{n_ki}, L_{n_kj})\ra \Hom_P(Pe_{\sum^{k-1}_{j=1}n_j+1}, L_{n_kj})\ra \Hom_P(Pe_{\sum^{k-1}_{j=1}n_j+i}, L_{kj})\ra 0.$$

Note that $\Hom_P(Pe_{\sum^{k-1}_{j=1}n_j+1}, L_{n_kj})$ and
$\Hom_P(Pe_{\sum^{k-1}_{j=1}n_j+i}, L_{kj})$ are both isomorphic to
$A_k/I_{kj}$ in $\Modcat{\mathbb{Z}}$. It follows from Lemma
\ref{144} that $\Hom_P(L_{n_ki}, L_{n_kj})=0$ for $2\leq i<j\leq
n_k, 1\leq k\leq m$.

Thus $\End_P((L_{n_12}\oplus\cdots\oplus
L_{n_1n_1})\oplus\cdots\oplus(L_{n_m1}\oplus\cdots\oplus
L_{n_mn_m})\oplus(Pe_{\sum^{m-1}_{l=1}n_l+1}\oplus\cdots\oplus
Pe_1))\cong$
$$\footnotesize\begin{pmatrix}
(L_{n_12}, L_{n_12})&0&0&\cdots&&&&0\\
*&\ddots&\ddots&\ddots&&&&\vdots\\
&&(L_{n_1n_1},L_{n_1,n_1})&0&&&&\\
&&\ddots&(L_{n_22},L_{n_22})&\ddots&&&\\
\vdots&&&&\ddots&&&\\
&&&&&(L_{n_2n_2},L_{n_2n_n})&&\\
&&&&&&\ddots&0\\
*&&&\cdots&&&*&\End_A(Ae_m\oplus\cdots\oplus Ae_1)
\end{pmatrix}.$$
By Lemma \ref{456} and Lemma \ref{PX}, we can get the
conclusion.$\square$

\begin{Koro}\label{cor51}
Let $\Lambda$ be as in Theorem \ref{Theo}. Then
$$\lfd (A)-1\leq \lfd(\Lambda)\leq n+\sum\limits_{i=2}^n\lfd
(A_i/I_i)+\lfd (A).$$
\end{Koro}
{\bf Proof.} It follows from Theorem \ref{Theo} and Lemma
\ref{PX}.$\square$

As an consequence of Theorem \ref{theo5}, we can get the following
corollary. By the corollary, we can get a class of algebras which
have finite finitistic dimension.

\begin{Koro}Let $\Lambda$ be as in Theorem \ref{Theo} and suppose that $A$ is an Artin
algebra. Then

$(1)$ If $\fd(A)<\infty$ and $\fd (A_j/I_{ji})<\infty$ for\, $2\leq
i\leq n_j, 1\leq j\leq m$, then $\fd (P)<\infty$.

$(2)$ If $\fd (P)<\infty$, then $\fd (A)<\infty$.
\end{Koro}

In \cite{O3}, K. Yamaura proved that any block extension of a basic
QF-algebra is a basic left Harada algebra. And for any basic left
Harada algebra $T$, there exists a basic QF-algebra $R$ such that
$T$ is isomorphic to an upper staircase factor algebra of a block
extension of $R$. By Theorem \ref{theo5}, we can get that the
finitistic dimension is finite for the block extension of a
QF-algebra. Thus, the finitistic dimension conjecture holds for the
class of left Harada algebras.

\begin{Koro}\label{cor512}
Suppose that $R$ is a QF-algebra and $P$ is the block extension of
$R$. Then
$$\fd (P)\leq \sum^m_{i=1}n_i-m<\infty.$$
\end{Koro}

{\bf Proof.} Note that $\fd (R)=0, \fd (A_j/\rad A_j)=0$ for $2\leq
i\leq n_j, 1\leq j\leq m$. $\square$

\begin{Prop}\label{lgd} Let $\Lambda$ as in Theorem \ref{Theo}. Then
$$\begin{array}{ll} max\{\lgd (A_i/I_i), \lgd (A), 2\leq i\leq n\}-1\leq
\lgd(\Lambda) \leq \displaystyle{\sum_{i=2,3,\cdots,n} \lgd
(A/I_i)}\\+\lgd (A)+n.
\end{array}$$
\end{Prop}

{\bf Proof.} By Theorem \ref{Theo}, we can get that the two rings
$\Lambda$ and $\Sigma$ are derived equivalent via a tilting module
whose projective dimension is less or equal 1. It follows from
\ref{comm} and Lemma \ref{456}.$\square$

\begin{Koro} Let $A$ be a ring with identity, $I_2, I_3, \cdots, I_n$
ideals of $A$. Set $$\Gamma=\small\begin{pmatrix}
  A & I_2& I_3& \cdots & \cdots & I_n\\
  A & A & I_3 & \cdots & \cdots & I_n\\
  A & I_2 & A &I_4 & \cdots & I_n\\
  A & I_2 & I_3 & A & \cdots & I_n\\
  \vdots & \vdots & \vdots && \ddots & \vdots \\
  A &I_2 & I_3 & \cdots & I_{n-1} & A
\end{pmatrix}$$

Then $\max\{ \lgd (A/I_i)-1, \lgd (A)-1, \pd(_AI_i), 2\leq i\leq
n\}\leq \lgd(\Gamma)\leq \max\{\lgd (A/I_i)+\pd(_AI_j)+3, \lgd
(A)+1, 2\leq i, j \leq n\}$.
\end{Koro}
{\bf Proof.} By Corollary 3.3(2) and Lemma \ref{456}, we can get the
conclusion. $\square$


\subsection{Tiled triangular rings}
Before we turn to the second topic, we recall the definition of
recollement, given by Beilinson, Bernstein and Deligne in their work
on perverse sheaves.

\begin{Def}\label{recollement}\cite{BBD}
Let $\mathcal {D}, \mathcal {D}'$ and $\mathcal {D}''$ be
triangulated categories. Then a recollement of $\mathcal {D}$
relative to $\mathcal {D}'$ and $\mathcal {D}''$, diagrammatically
expressed by
$$\xymatrix{\mathcal {D}'
\ar[r]|{} &\mathcal {D}\ar@<1.2ex>[l]^{}\ar@<-1.2ex>[l]_{}\ar[r]|{}
&\mathcal {D}''\ar@<1.2ex>[l]^{}\ar@<-1.2ex>[l]_{}} $$ is given by
six exact functors
$$i_*=i_!: \mathcal {D}'\ra \mathcal {D}, j^*=j^!: \mathcal {D}\ra\mathcal {D}'',
i^*, i^!: \mathcal {D}\ra\mathcal {D}', j_!, j_*: \mathcal {D}''\ra
\mathcal {D},$$ which satisfy the following four conditions:

$(R1)$ $(i^*, i_*=i_!, i^!)$ and $(j_!, j^*=j^!, j_*)$ are adjoint
triples, i.e., $i^*$ is left adjoint to $i_*$ which is left adjoint
to $i^!$ etc.,

$(R2)$ $i^!j_*=0$,

$(R3)$ $i_*, j_!$ and $j_*$ are full embeddings,

$(R4)$ any object $X$ in $\mathcal {D}$ determines distinguished
triangles
$$i_!i^!X\ra X\ra j_*j^*X\ra \Sigma i_!i^!X \mbox{ and }
j_!j^!X\ra X\ra i_*i^*X\ra\Sigma j_!j^!X$$ where the morphisms
$i_!i^!X\ra X$ etc. are the adjunction morphisms.
\end{Def}

Using the notion of recollement, Happel proved the following result.
The next lemma is useful to provide a class of algebras which have
finite finitistic dimension.

\begin{Lem}\cite{H1}\label{dh} Let $A$ be a finite-dimensional algebra and
assume that $\D^b(\modcat{A})$ has a recollement  relative to
$\D^b(\modcat{A'})$ and $\D^b(\modcat{A''})$ for some finite
dimensional algebras $A'$, $A''$. Then $\fd(A)<\infty$ if and only
if $\fd (A')<\infty$ and $\fd(A'')<\infty$.
\end{Lem}

The following lemma, showing how to construct a recollement, is
useful in our proof.








\begin{Lem}\cite{NS}\label{NS}Let $A, B$ and $C$ be algebras. The
following assertions are equivalent:

$(1)$ $\D^-(A)$ is a recollement of $\D^-(C)$ and $\D^-(B)$.

$(2)$ There are two objects $P, Q\in \D^-(A)$ satisfying the
following properties:

\quad$(a)$ There are isomorphism of algebras $C\cong
\Hom_{\D(A)}(P,P)$
 and $B\cong \Hom_{\D(A)}(Q,Q)$.

\quad$(b)$ $P$ is exceptional and isomorphic in $\D(A)$ to a bounded
complex of finitely generated projective $A$-modules.

\quad$(c)$ For every set $\Lambda$ and every non-zero integer $i$ we
have $\Hom_{\D(A)}(Q, Q^{(\Lambda)}[i])=0$, the canonical
isomorphism $\Hom_{\D(A)}(Q,Q)^{(\Lambda)}\ra
\Hom_{\D(A)}(Q,Q^{(\Lambda)})$ is an isomorphism, and $Q$ is
isomorphic in $\D(A)$ to a bounded complex of projective
$A$-modules.

\quad$(d)$ $\Hom_{\D(A)}(P,Q[i])=0$ for all $i\in \mathbb{Z}$.

\quad$(e)$ $P\oplus Q$ generates $\D(A)$.
\end{Lem}






Now, we turn to consider ``tiled triangular ring,'' i.e., rings of the form\\
$$\Delta=\small
\begin{pmatrix}
A&I_{12}&\cdots&I_{1n}\\
A&A&&\vdots\\
\vdots&&\ddots&I_{n-1,n}\\
A&\cdots&\cdots&A
 \end{pmatrix}$$
for $I_{ij}$ ideals of $A$.

Now, let us prove the last result in this paper.

\begin{Prop}\label{fd} Set
$$\Phi =
 \small\begin{pmatrix}
A & I_{1,2} &I_{1,3} & \cdots & I_{1,n-1} & I_{1,n}\\
A & A & I_{2,3} & \cdots & I_{2,n-1} & I_{2,n}\\
A & A& A& \cdots & I_{3,n-1} & I_{3,n}\\
A & A & A& A & \cdots & I_{4,n}\\
\vdots & \vdots &\ddots &\ddots  & \ddots& \vdots\\
A & A & \cdots & \cdots & A & A
 \end{pmatrix}$$

Suppose that $\Phi$ is an Artin algebra, $I_{ij}$ are ideals of $A$
for $1\leq i<j\leq n$,
$\pd(_{A/I_{i,i+1}}I_{i+1,j+1}/I_{i,j+1})<\infty$ for $1\leq i<j\leq
n-1$ and $\fd(A/I_{i,i+1})<\infty, \fd(A)<\infty$ for $1\leq i\leq
n-1$. Then $\fd(\Phi)<\infty$.
\end{Prop}
{\bf Proof.} Set $\Gamma=M_n(A)$.
The exact sequence
$$0\lra \Phi \stackrel{\lambda}\lra \Gamma\stackrel{\pi}\lra L\lra 0 $$
in $\Phi$-Mod is an $\add(\Gamma)$-split sequence. By the method
which is similar to the ones in Theorem \ref{Theo}, we can prove
that the two rings
$$\Phi=\small\begin{pmatrix}A&I_{12}&I_{13}&\cdots&I_{1,n}\\
A&A&I_{23}&&\vdots\\
\vdots&\vdots&\ddots&\ddots&\vdots\\
\vdots&\vdots&&\ddots&I_{n-1,n}\\
A&A&A&A&A
\end{pmatrix}\quad\mbox{ and }\quad
\Sigma_1=\small\begin{pmatrix}A/I_{12}&I_{23}/I_{13}&\cdots&I_{2n}/I_{1n}&0\\
A/I_{12}&A/I_{13}&&&\vdots\\
\vdots&\vdots&&\ddots&\vdots\\
A/I_{12}&A/I_{13}&\cdots&A/I_{1n}&0\\
A/I_{12}&A/I_{13}&\cdots&A/I_{1n}&A
\end{pmatrix}$$ are derived equivalent.

Set
$$\Phi_{n-1}=\small\begin{pmatrix}A/I_{12}&I_{23}/I_{13}&\cdots&I_{2n}/I_{1n}\\
A/I_{12}&A/I_{13}&\ddots&\vdots\\
\vdots&\vdots&\ddots&I_{n-1,n}/I_{1n}\\
A/I_{12}&A/I_{13}&\cdots&A/I_{1n}
\end{pmatrix}$$

For simplicity, we denote $\Phi_{n}$ by
$$\Gamma=\small\begin{pmatrix}A_1&J_{12}&J_{13}&\cdots&J_{1n}\\
\vdots&A_2&J_{23}&&\vdots\\
\vdots&\vdots&A_3&\ddots&\vdots\\
\vdots&\vdots&\vdots&\ddots&J_{n-1,n}\\
A_1&A_2&A_3&\cdots&A_n.
\end{pmatrix}$$

{\bf Claim:} Suppose that $\pd(_{A_1}J_{1k})<\infty, 2\leq k\leq n$
and $\pd(_{{A_i/J_{i-1,i}}}J_{ij}/J_{i-1,j})<\infty, 2\leq i<j\leq
n-1$. Then $\fd(\Gamma) <\infty$ if and only if $\fd(A_1)\leq\infty$
and $\fd(A_i/J_{i-1,i})<\infty$, for $2\leq i\leq n$.

{\bf Proof of Claim:} Let $e$ be an idempotent of $\Gamma$ which has
$1$ in the $(1,1)$-th position and zeros elsewhere. By easy
computation,
$$\Gamma
e\Gamma=\small\begin{pmatrix}A_1&J_{12}&J_{13}&\cdots&J_{1n}\\
A_1&J_{12}&J_{13}&\cdots&J_{1n}\\
\vdots&&\vdots&&\vdots\\
A_1&J_{12}&J_{13}&\cdots&J_{1n}
\end{pmatrix}.$$

Since $\Gamma e\Gamma$ is projective as right $\Gamma$-module, we
have $\Tor^{\Gamma}_i(\Gamma/\Gamma e\Gamma, \Gamma/\Gamma
e\Gamma)=0$ for $i>0$. Then $\lambda:{} \Gamma\ra \Gamma/\Gamma
e\Gamma$ is a homological ring homomorphism. Then there is
recollement:
$$\xymatrix{\D(\Mod-\Gamma/{\Gamma
e\Gamma})\ar[r]|(.6){i_*} &\D(\Mod-\Gamma)\ar@<1.2ex>[l]^{{\bf
R}\Hom_{\Gamma}(\Gamma/\Gamma
e\Gamma,-)}\ar@<-1.2ex>[l]_{-\otimes^{\bf L}_\Gamma \Gamma/{\Gamma
e\Gamma}}\ar[r]|{j^!} &{\rm Tria}_{\D(\Gamma)}(\Gamma
e\Gamma)\ar@<1.2ex>[l]^{j_*}\ar@<-1.2ex>[l]_{-\otimes^{\bf L}_\Gamma
\Gamma e\Gamma}}$$ where ${\rm Tria}_{\D(\Gamma)}(\Gamma e\Gamma)$
is the smallest full triangulated subcategory of $\D(\Gamma)$
containing $\Gamma e\Gamma$ and closed under small coproducts, $i_*$
is the inclusion functor, $-\otimes^{\bf L}_{\Gamma}\Gamma/\Gamma
e\Gamma $, is the left derived functor of
$-\otimes_{\Gamma}\Gamma/\Gamma e\Gamma$, $-\otimes^{\bf
L}_{\Gamma}\Gamma e\Gamma$ is the left derived functor of
$-\otimes_{\Gamma}\Gamma e\Gamma$ and ${\bf
R}\Hom_{\Gamma}(\Gamma/\Gamma e\Gamma,-)$ is the right derived
functor of $\Hom_{\Gamma}(\Gamma/\Gamma e\Gamma,-)$.

Note that $\Gamma/\Gamma e\Gamma$ has finite projective dimension as
right $\Gamma$-module. Then, by Lemma \ref{NS} and \cite[Corollary 3
and Example 6]{NS}, there is a recollement:
$$\xymatrix{\D^-(\Mod-\Gamma/{\Gamma
e\Gamma})\ar[r]|(.6){i_*} &\D^-(\Mod-\Gamma)\ar@<1.2ex>[l]^{{\bf
R}\Hom_{\Gamma}(\Gamma/\Gamma
e\Gamma,-)}\ar@<-1.2ex>[l]_{-\otimes^{\bf L}_\Gamma \Gamma/{\Gamma
e\Gamma}}\ar[r]|{}
&\D^-(C)\ar@<1.2ex>[l]^{}\ar@<-1.2ex>[l]_{-\otimes^{\bf L}_\Gamma
\Gamma e\Gamma}}$$ where $C$ is the dg algebra $\mathcal
{C}_{dg}\Gamma(i\Gamma e\Gamma, i\Gamma e\Gamma)$ and $i\Gamma
e\Gamma$ is an injective resolution of the right $\Gamma$-module
$\Gamma e\Gamma$. Note that $\Gamma e\Gamma$ is isomorphic to
$(e\Gamma)^n$ as right $\Gamma$-module. By \cite[Theorem 9.2]{K1},
there is a triangle equivalence between $\D^-(C)$ and
$\D^-(\Mod-H^0(C))=\D^-(\Mod-e\Gamma e)$. Hence there is a
recollement:
$$\xymatrix{\D^-(\Mod-\Gamma/{\Gamma
e\Gamma})\ar[r]|(.6){i_*} &\D^-(\Mod-\Gamma)\ar@<1.2ex>[l]^{{\bf
R}\Hom_{\Gamma}(\Gamma/\Gamma
e\Gamma,-)}\ar@<-1.2ex>[l]_{-\otimes^{\bf L}_\Gamma \Gamma/{\Gamma
e\Gamma}}\ar[r]|{j^!} &\D^-(\Mod-e\Gamma
e)\ar@<1.2ex>[l]^{j_*}\ar@<-1.2ex>[l]_{-\otimes^{\bf L}_{e\Gamma
e}e\Gamma}}$$ where $i_*$ is inclusion functor, $j^!={\bf
R}\Hom_{\Gamma}(e\Gamma, -)$ is the right derived functor of
$\Hom_{\Gamma}(e\Gamma,-)$, $j_*={\bf R}\Hom_{e\Gamma e}(\Gamma
e,-)$ is the right derived functor of $\Hom_{e\Gamma e}(\Gamma
e,-)$, $-\otimes^{\bf L}_{\Gamma}\Gamma/\Gamma e\Gamma$ is the left
derived functor of $-\otimes_{\Gamma}\Gamma/\Gamma e\Gamma$,
$-\otimes^{L}_{e\Gamma e}e\Gamma$ is the left derived functor of
$-\otimes_{e\Gamma e}e\Gamma$. Since $\pd(_{A_1}J_{1k})<\infty,
2\leq k\leq n$, we have that $\Gamma e\Gamma$ have finite projective
dimension as left $\Gamma$-module. Thus, the functors $-\otimes^{\bf
L}_{\Gamma}\Gamma/\Gamma e\Gamma$ and $-\otimes^{\bf L}_{e\Gamma e}
e\Gamma$ send complexes of bounded homology to complexes of bounded
homology. Note that $e\Gamma$ and $\Gamma/\Gamma e\Gamma$ have
finite projective dimension as right $\Gamma$-module. Then the
functors ${\bf R}\Hom_{\Gamma}(e\Gamma, -)$ and ${\bf
R}\Hom_{\Gamma}(\Gamma/\Gamma e\Gamma,-)$ restrict to the functors
$\D^b(\Gamma-\Mod)\ra \D^b(e\Gamma e-\Mod)$ and
$\D^b(\Gamma-\Mod)\ra \D^b(\Gamma/\Gamma e\Gamma-\Mod)$
respectively. Then, we can get a recollement:
$$\xymatrix{\D^b(\Mod-\Gamma/{\Gamma
e\Gamma})\ar[r]|(.6){i_*} &\D^b(\Mod-\Gamma)\ar@<1.2ex>[l]^{{\bf
R}\Hom_{\Gamma}(\Gamma/\Gamma
e\Gamma,-)}\ar@<-1.2ex>[l]_{-\otimes^{\bf L}_\Gamma \Gamma/{\Gamma
e\Gamma}}\ar[r]|{j^!} &\D^b(e\Gamma
e-\Mod)\ar@<1.2ex>[l]^{j_*}\ar@<-1.2ex>[l]_{-\otimes^{\bf
L}_{e\Gamma e}e\Gamma}}.\quad\quad(*)$$

Since $\Gamma$ is Artin algebra and the functors, appearing in
$(*)$, take finitely generated modules to finite generated modules,
it follows that the following diagram is a recollement.
$$\xymatrix{\D^b(\Gamma/\Gamma
e\Gamma)\ar[r]|{}
&\D^b(\Gamma)\ar@<1.2ex>[l]^{}\ar@<-1.2ex>[l]_{}\ar[r]|{}
&\D^b(e\Gamma e)\ar@<1.2ex>[l]^{}\ar@<-1.2ex>[l]_{}}.$$

By Lemma \ref{dh}, $\fd(\Gamma)<\infty$ if and only if $\fd
(A_1)<\infty$ and $\fd(\Gamma/{\Gamma e\Gamma})<\infty$. Let
$\Gamma_{n-1}$ denote $\Gamma/{\Gamma e\Gamma}$. By similar
discussion, we can get $\fd(\Gamma)<\infty$ if and only if $\fd
(A_1)<\infty$ and $\fd(A_i/J_{i-1,i})<\infty$ for $2\leq i\leq
n$.$\square$

By Claim, $\fd(\Phi_{n-1})<\infty$ if and only if $\fd
(A/I_{i,i+1})<\infty$ for $1\leq i\leq n-1$. By Lemma \ref{456},
$\fd(A)\leq \fd(\Sigma_1) \leq \fd(A)+\fd(\Phi_{n-1})+1$. By Claim,
$\fd(A/I_{i,i+1})<\infty, \fd(A)<\infty, 1\leq i\leq n-1$ imply
$\fd(\Sigma_1)<\infty$. By Lemma \ref{PX}, we can get that
$\fd(\Sigma_1)<\infty$ implies $\fd(\Phi)<\infty$.$\square$

\medskip

The following is a typical case of Theorem \ref{fd}.

\begin{Koro}\label{last} Set
$$\Phi =
 \small\begin{pmatrix}
A & \rad A &I_{1,3} &\cdots && I_{1,n}\\
A & A & \rad A &\ddots & &\vdots\\
A & A& A&\ddots & &\vdots\\
\vdots & \vdots &\ddots & & \ddots& I_{n-2,n}\\
\vdots&\vdots&&&\ddots&\rad A\\
A & A & \cdots &  & A & A
 \end{pmatrix}$$

Suppose that $\Phi$ is Artin algebra, $I_{ij}$ are ideals of $A$ for
$i,j=1,2,\cdots,n$. If $\fd(A)<\infty$, then $\fd(\Phi)<\infty$.
\end{Koro}

{\bf Proof.} Let $I_{i,i+1}=\rad A$ for $i=1,2,\cdots,n-1$.$\square$


\section{Examples} In this section, we display several examples to illustrate our theorem.

\noindent{\bf Example 1.} Let $A$ be a ring with  ideal $I$. Let
$\Gamma$ be the ring
$$\Gamma=\begin{pmatrix}
  A&I&I^2&I^3\\
  A&A&I^2&I^3\\
  A&I&A&I^3\\
  A&I&I&A
\end{pmatrix} $$

By Proposition \ref{lgd}, we can get $max\{0,\lfd(A)-1\}\leq \lfd\Gamma
\leq \lfd(A)+\lfd(A/I)+\lfd(A/I^2)+\lfd(A/I^3)+4$.

\medskip

\noindent{\bf Example 2.} Let $A=k[x]/(x^n)$ for $n\geq 1$, and
$I=\rad(A)$, Let
$$\Lambda=\begin{pmatrix}
A&I&I^2&I^3\\
A&A&I&I^3\\
A&A&A&I\\
A&A&A&A\\
\end{pmatrix}$$


Since $k[x]/(x^n)$ is representation-finite, the finitistic
dimension of $k[x]/(x^n)$ is finite. By Corollary \ref{last},
$\fd(\Lambda)<\infty$.
\medskip

\noindent{\bf Example 3.} Let $A$ be a $k$-algebra given by the
following quiver.
$$
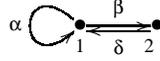
\xymatrix{
  *{\bullet}\ar[r]^{\beta} \ar@(ul,dl)_{\alpha}& *{\bullet}\ar@<2pt>[l]^(0.1){2}^(1){1}^{\delta}  }$$
with relations $\{\alpha^3=\beta\delta, \alpha\beta=0,
\delta\alpha=0\}$.

Then $A$ can be represented as the following matrix form.

$$A=\begin{pmatrix}
  k[\alpha]/(\alpha)^4&k\beta\\
  k\delta&k[\delta\beta]/(\delta\beta)^2
\end{pmatrix}.$$

Suppose that $P(3, 2)$ is the block extension of $A$. $P(3, 2)=$
$$\footnotesize
\begin{pmatrix}
k[\alpha]/(\alpha)^4&k[\alpha]/(\alpha)^4&k[\alpha]/(\alpha)^4&k\beta&k\beta\\
(\alpha)/(\alpha)^4&k[\alpha]/(\alpha)^4&k[\alpha]/(\alpha)^4&k\beta&k\beta\\
(\alpha)/(\alpha)^4&(\alpha)/(\alpha)^4&k[\alpha]/(\alpha)^4&k\beta&k\beta\\
k\delta&k\delta&k\delta&k[\delta\beta]/(\delta\beta)^2&k[\delta\beta]/(\delta\beta)^2\\
k\delta&k\delta&k\delta&(\delta\beta)/(\delta\beta)^2&k[\delta\beta]/(\delta\beta)^2\\
\end{pmatrix}$$

$P(3,2)$ can be described by given by the quiver
$$
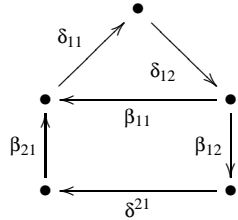
\xymatrix{
  & \bullet\ar[rd]_{\delta_{12}}  &  \\
  \bullet  \ar[ur]^{\delta_{11}} &  & \bullet \ar[d]_{\beta_{12}}\ar[ll]^{\beta_{11}} \\
  \bullet \ar[u]^{\beta_{21}} &  & \bullet\ar[ll]^{\delta^{21}}  }$$
with relations
$\{e(\alpha^3)=\beta_{11}\delta_{11}\delta_{12}\beta_{11}\delta_{11}\delta_{12}\beta_{11}=
\beta_{12}\delta_{21}\beta_{21}=e(\beta\delta),
e(\alpha\beta)=\beta_{11}\delta_{11}\delta_{12}\beta_{12}=0,
e(\delta\alpha)=\beta_{21}\delta_{11}\delta_{12}\beta_{11}=0\}$,
where $e$ is the extension map defined in \cite{KY}.
By Corollary \ref{cor512}, we have $\fd P(3, 2)\leq 3$.

\medskip {\bf Acknowledgement.}
I want to express my gratitude to my
supervisor professor Changchang Xi for encouragement and useful suggestions.

\medskip \footnotesize{
}


\begin{thebibliography}{99}


\bibitem{ARS}{{\sc M. Auslander and I. Reiten and S. O. Smal{\o},} {\it Representation Theory of Artin Algebras.} Cambridge University Press. 1995.}

\bibitem{BBD}{{\sc A. A. Beilinson, J. Bernstein and P. Deligne,} Faisceaux pervers, in Analyse et topologie sur les espaces
singuliers. Ast\'{e}risque 100 (1982), 1-172.}

\bibitem{YPC1}{{\sc Y. P. Chen} Derived equivalences in n-angulated categries. arXiv: 1107.2985.}

\bibitem{YPC2}{{\sc Y. P. Chen} Constructions of derived equivalences between subrings. Preprint.}

\bibitem{C}{{\sc K. N. Cowley,} On the global dimension of matrix subrings.
Comm. Algebra 24 (1996), 857-871.}

\bibitem{DS}{{\sc D. Dugger and B. Shipley,} K-theory and derived
equivalences. Duke Math. J. 124 (3)(2004), 587-617.}



\bibitem{FGR}{{\sc R. M. Fossum, P. A. Griffith and I.Reiten,}
{\it Trivial Extensions of Abelian Categories.} Springer-Verlag,
Berlin, 1975.}




\bibitem{H}{{\sc D. Happel,} {\it Triangulated Categories in the Representation Theory of Finite Dimensional
Algebras.} Cambridge Univ. Press, Cambridge, 1988.}
\bibitem{H1}{{\sc D. Happel,} Reduction techniques for homological conjecture.
Tsukuba J. Math 17 (1993), 115-130.}


\bibitem{HX1}{{\sc W. Hu and C. C. Xi,} $\mathcal {D}$-split sequences and derived equivalences. Adv. Math. 227 (2011), 292-318.}


\bibitem{K}{{\sc Y. Kato,} On derived equivalent coherent rings.
Comm. Algebra 30 (2002), 4437-4454.}
\bibitem{K1}{{\sc B. Keller,} Deriving DG categories. Ann. Sci. \'{E}cole Norm. Sup. (4)27(1994), no. 1, 63-102.}

\bibitem{KK}{{\sc E. Kirkman and J. Kuzmanovich,} Matrix subrings having finite global dimension.
J. Algebra 109 (1987), 74-92.}
\bibitem{KKS}{{\sc E. Kirkman, J. Kuzmanovich and L. Small,} Finitistic dimension of noetherian rings.
J. Algebra 147 (1992), 350-364.}






\bibitem{MB}{{\sc J. C. McConnel and J. C. Bobson,} {\it Noncommutative Noetherian Rings}.
 Wiley, Chichester, England, 1987.}



\bibitem{NS}{{\sc P. Nicol\'{a}s and M. Saor\'{\i}n,} Lifting and restricting recollement data. Appl. categor. struct.,
DOI 10. 1007/s10485-009-9198-z.}



\bibitem{O3}{{\sc K. Oshiro,} On Harada rings III. Math. J. Okayama Univ. 32 (1990), 111-118.}



\bibitem{PX}{{\sc S. Y. Pan and C. C. Xi,} Finiteness of finitistic dimension is invariant of derived equivalences. J. Algebra 322 (2009), 21-24.}


\bibitem{R1}{{\sc J. Rickard,} Morita theory for derived categories. J. London Math. Soc. 39(1989), 436-456.}
\bibitem{R2}{{\sc J. Rickard,} Derived categories and stable equivalences. J. Pure Appl. Algebra 64(1989), 303-317. }





\bibitem{XX}{{\sc C. C. Xi and D. M. Xu,} The finitistic dimension conjecture and relatively
projective modules. Preprint is available at: http:
//math.bnu.edu.cn/$^{\sim}$ccxi/.}

\bibitem{KY}{{\sc K. Yamaura,} Quivers with relations of Harada algebras. Proc. Amer. Math. Soc. 138 (2010), 47-59.}

\end{thebibliography}
\end{document}